\documentclass[12pt,a4paper,fleqn]{article}

\usepackage{a4wide}
\usepackage{amsmath,amssymb,amsthm}
\theoremstyle{plain}
\newtheorem{theorem}{Theorem}
\newtheorem*{theorem*}{Theorem}
\newtheorem{corollary}[theorem]{Corollary}

\newtheorem{lemma}[theorem]{Lemma}
\newtheorem{proposition}[theorem]{Proposition}

\theoremstyle{definition}
\newtheorem{definition}[theorem]{Definition}

\newtheorem*{note*}{Note}
\newtheorem*{notes*}{Notes}
\newtheorem*{null*}{}
\newtheorem{remark}[theorem]{Remark}
\usepackage[no-weekday]{eukdate}
\usepackage[T1]{fontenc}
\usepackage[pdftex]{graphics}
\usepackage[pdftex,pagebackref=true,colorlinks,linkcolor=blue,citecolor=blue]{hyperref}
\renewcommand*{\backref}[1]{}
\renewcommand*{\backrefalt}[4]{[{%
    \ifcase #1 Not cited.%
    \or Cited on page~#2.%
    \else Cited on pages #2.%
    \fi%
  }]}
\usepackage{ifthen}
\usepackage{mathpartir}
\usepackage[numbers]{natbib}
\usepackage{stmaryrd}
\usepackage{url}
\usepackage[all,pdf]{xy}
\CompileMatrices 

\usepackage[tt=false, type1=true]{libertine}
\usepackage[libertine]{newtxmath}
\usepackage[sans]{dsfont}
\usepackage[mathscr]{eucal}

\newcommand{\aczel}{\mathrel{\lhd}}
\newcommand{\bimp}{\mathbin{\leftrightarrow}}
\newcommand{\Bimp}{\Leftrightarrow}

\newcommand{\comp}{\circ}
\newcommand{\conj}{\wedge}
\newcommand{\Conj}{\mathrel{\&}}

\newcommand{\defeq}{\triangleq}
\newcommand{\disj}{\vee}

\renewcommand{\d}{\mathrm{d}}
\newcommand{\down}{{\downarrow}}
\newcommand{\E}{\mathscr{E}}

\newcommand{\epi}{\twoheadrightarrow}
\newcommand{\Fin}{\mathrm{Fin}}
\renewcommand{\H}{\mathbf{H}}
\newcommand{\Heyt}{\mathbf{Heyt}}

\newcommand{\imp}{\mathbin{\shortrightarrow}}
\newcommand{\Imp}{\Rightarrow}
\newcommand{\iso}{\cong}

\newcommand{\join}{\vee}
\renewcommand{\Join}{{\vee}}
\newcommand{\km}{\partial}
\newcommand{\KM}{\nabla}
\newcommand{\meet}{\wedge}
\newcommand{\Meet}{{\wedge}}
\newcommand{\mono}{\rightarrowtail}
\newcommand{\morphism}{\rightarrow}
\newcommand{\NNO}{\mathbb{N}}

\newcommand{\p}{\mathrm{p}}

\newcommand{\scn}{\mathop{\mathbf{scn}}}
\newcommand{\Set}{\mathbf{Set}}

\newcommand{\size}{\mathop{size}}

\newcommand{\up}{{\uparrow}}

\title{A Note on the Topos-Inadmissibility\\
  of Free Heyting Algebras}

\author{Andrew Pitts\\
  University of Cambridge}

\date{\today{}}

\begin{document}

\maketitle

\begin{abstract}
  Yiqi Xu and Linguan Ye have recently settled the long-standing open
  question about which Heyting algebras are \emph{topos-admissible},
  that is, occur as the truth values of some topos. They prove that
  free Heyting algebras with two or more generators are not topos
  admissible. Their proof has been simplified by Zoltan Kocsis and by
  Sam van~Gool. Here I make some further simplifications, emphasising
  algebraic methods over Kripke semantics and in particular the use of
  Sierpi{\'n}ski cones of Heyting algebras.
\end{abstract}

\section{Introduction}

If $\E$ is a (locally small) topos~\cite{JohnstonePT:skeett}, then its
subobject classifier $\Omega_\E$ is a Heyting algebra object in $\E$
and hence the set of its truth values $1\morphism\Omega_\E$ (in
bijection with the subobjects of the terminal object $1$ in $\E$) is a
Heyting algebra in the category of sets and functions. Following
\citet{AbbadiniM:topefh}, call a Heyting algebra
\emph{topos-admissible} if it is isomorphic to one of this form and
\emph{topos-inadmissible} if it is not.  Complete Heyting algebras are
topos-admissible via toposes of sheaves; Boolean algebras are
topos-admissible~\citep[Exercise~9.11]{JohnstonePT:topt};
and recently \citet{AbbadiniM:topefh} have shown, more generally, that
all {étale}-finite Heyting algebras are topos-admissible. The question
of whether \emph{all} Heyting algebras are in fact topos-admissible is
probably as old as Lawvere and Tierney's notion of elementary topos
(1969--70). I had my attention drawn to the question by a seminar on
Heyting algebras given by Freyd while he was visiting Cambridge in
October 1980.

\citet{XuY:faihot} have settled in the negative this long-standing
question, by proving that \emph{the free Heyting algebra on two
  generators is topos-inadmissible}. Write $\H(X)$ for the free
Heyting algebra on a set $X$ of generators. In broad terms, Xu and
Ye's proof strategy is to demonstrate that $\H(2)$ fails to have
certain completeness properties that it would have to have were it
topos-admissible. Ironically this is no different from
the strategy I adopted in the 1980s when trying to show that
$\H(\omega)$ is topos-inadmissible by proving that the inclusion
$\H(\omega)\mono\H(\omega+1)$ does not have left and right adjoints
(it is not too hard to see that such adjoints would have to exist were
$\H(\omega)$ to be topos-admissible). Contrary to expectations, those
adjoints do exist~\citep{PittsAM:intsoq}.  Although that means this
particular completeness property does not serve the purpose, the
consolation prize was an interesting new property of intuitionistic
and other logics, uniform interpolation~\citep{vanGoolS:unii}.

During the early part of this century the expectation was that the
question about topos admissibility has a positive answer. The reason
for this was an incomplete preprint by the late~Dimitri
Pataraia~\citep{PataraiaD:higoqt} that circulated among interested
parties. Its abstract claims to resolve the question via the notion of
``higher order cylindric Heyting algebra'' (\emph{hocha}). This gives
a cylindric algebra treatment of intuitionistic higher-order logic, in
contrast to the treatment in terms of
triposes~\citep{PittsAM:trit,PittsAM:tritr} founded upon Lawvere's
category-theoretic notion of
hyperdoctrine~\citep{LawvereFW:equhcs}. Hochas correspond to minimal
toposes~\citep{JohnstonePT:hocmt} and in particular, given a hocha one
can construct a corresponding topos whose Heyting algebra of truth
values is isomorphic to the one given by the hocha. Thus it remains to
construct a hocha with a given Heyting algebra of
truth-values. Pataraia [private communication] intended to do this via
the particular construction of free Heyting algebras due to
\citet{GhilardiS:frehab}, but the details are missing from
\citep{PataraiaD:higoqt} and now we know that this approach cannot
work in general.

In 2025 the talk~\citep{PittsAM:heyaho} encouraged renewed attention
to this question of topos admissibility, resulting in work such
as \citep{AbbadiniM:topefh} and culminating in the negative solution by
\citet{XuY:faihot}. Their proof of topos-inadmissibility of $\H(2)$
makes use of Bellissima's construction of finitely generated free
Heyting algebras using Kripke
models~\citep{BellissimaF:fingfh,DarniereL:belcfg}
\citep[Remark~4.17]{YeL:coecph}.  The proof has been simplified by
\citet{KocsisZ:twoapf} and by \citet{vanGoolS:gentvt}, also using
Kripke semantics. In particular Kocsis shows that an incompleteness
property of $\H(2)$ that does serve the purpose has to do with the
following subset of a Heyting algebra $H$ equipped with a
distinguished element $a$:
\[
  \KM(H,a) \defeq \{h \disj(h\imp a) \mid h \in H\} 
\]
\citet{KocsisZ:twoapf} cites as motivation for considering $\KM(H,a)$
the work of \citet{SimmonsH:algvcb} on a localic version of the
Cantor--Bendixson coderivative and the recent work of
\citet{FereeH:pitims} on uniform interpolation for the KM modal
logic. When $H=\H(2)$ and $a$ is one of the two generators, one can
show that (a) \emph{$\KM(H,a)$ would have to have an infimum were
  $\H(2)$ to be topos-admissible}, but that (b) \emph{it does not have an
  infimum}. Here I give proofs of (a) and (b) emphasising algebraic
methods over Kripke semantics and in particular the use of
Sierpi{\'n}ski cones of Heyting algebras.

\section{Preliminaries on Heyting Algebras}

In this note $\Heyt$ is the category of Heyting algebras and their
morphisms~\citep[Chapter~IX]{BalbesR:disl}. The partial order, top,
binary meet, bottom, binary join, Heyting implication and Heyting
bi-implication operations of a Heyting algebra are denoted by $\leq$,
$\top$, $\conj$, $\bot$, $\disj$, $\imp$ and $\bimp$ respectively.  I
write $h < k$ to mean $h\leq k \Conj h\not= k$. Meets and joins of
subsets $S$ of a Heyting algebra are denoted by $\Meet S$ and
$\Join S$ when they exist---for example, when $S$ is finite, or when
the Heyting algebra is complete.

\subsection{Free Heyting algebras}
\label{sec:freha}

Let $\H(X)$ denote the \emph{free Heyting algebra} on a set $X$. It
comes equipped with a function
\begin{equation}
  \label{eq:18}
  \eta_X : X \morphism \H(X)
\end{equation}
that together are determined up to isomorphism by the universal
property that for all functions $f:X \morphism H$ from $X$ to a
Heyting algebra $H$, there is a unique morphism of Heyting algebras
$\bar{f}:\H(X)\morphism H$ with $\bar{f}\comp \eta_X = f$. Since there
is an injection from $X$ to some Heyting algebra (for example the
powerset of $X$ ordered by inclusion), $\eta_X$ is an
injection. Sometimes I identify $X$ with its image under $\eta_X$ and
just write $x \in \H(X)$ when $x \in X$.

Recall that the initial Heyting algebra $\H(\emptyset)$ is the
two-element algebra $2 = \{\bot < \top\}$ and that the free Heyting
algebra on one generator $\H(1)$ is the Rieger--Nishimura
lattice~\citep{RiegerL:lattbp,NishimuraI:forvip}, which is an
infinite, but complete Heyting algebra.

\subsection{Terms}
\label{sec:ter}

Free Heyting algebras can be presented as quotients of the set
$\Sigma(X)$ of algebraic terms over a set $X$ of generators for a
signature $\Sigma$, for a number of different choices of signature
$\Sigma$. It is convenient (at least for the proof of
Proposition~\ref{prop:cof} in the Appendix) to use a signature with
operations for all finite meets and joins rather than just for nullary
and binary ones. So inductively define $\Sigma(X)$ as follows:
\begin{itemize}
\item each $x\in X$ is in $\Sigma(X)$
  
\item for each finite subset $S\subseteq \Sigma(X)$ we have  a
  conjunction term
  $\Meet S$ and a disjunction term $\Join S$ in $\Sigma(X)$
  
\item for each $\varphi$ and $\psi$ in $\Sigma(X)$ there is an
  implication term $\varphi \imp \psi$ in $\Sigma(X)$.
\end{itemize}
Then $\H(X)$ is the quotient of $\Sigma(X)$ by the equivalence
relation associated with the usual pre-order on $\Sigma(X)$ for
intuitionistic entailment, here just denoted by $\leq_X$. Write
$[\varphi]_X\in \H(X)$ for the equivalence class of
$\varphi \in \Sigma(X)$. Abbreviate $\leq_X$ and $[\varphi]_X$ by
$\leq$ and $[\varphi]$ in case $X$ is clear from the context.

The operations on $\Sigma(X)$ for falsity ($\bot$), binary disjunction
($\varphi\disj \psi$), truth ($\top$), binary conjunction
($\varphi\conj \psi$) and intuitionistic bi-implication
($\varphi\bimp\psi$) are treated as abbreviations:
$\bot \defeq \Join\emptyset$,
$\varphi\disj \psi \defeq \Join\{\varphi,\psi\}$,
$\top \defeq \Meet\emptyset$,
$\varphi\conj \psi \defeq \Meet \{\varphi,\psi\}$,
$\varphi\bimp\psi \defeq \Meet\{\varphi\imp\psi, \psi\imp\varphi\}$.

\smallskip

The following property goes back to
\citet{McKinseyJCC:cloeca,McKinseyJCC:thescl}.
\begin{proposition}[Finite model property]
  \label{prop:fmp}
  For each free Heyting algebra $H$, the collection of surjective
  morphisms\footnote{Recall that these are exactly the epimorphisms in
    $\Heyt$.}  $H\epi K$ in $\Heyt$ with $K$ finite, is jointly monic.
\end{proposition}
\begin{proof}
  Suppose $H= \H(X)$ for some set $X$. It suffices to construct for
  each $\varphi\in\Sigma(X)$ a surjective morphism
  $f_\varphi: \H(X)\morphism K_\varphi$ to a finite Heyting algebra
  with the property that
  ${f_\varphi[\varphi] = \top} \Imp {[\varphi]=\top}$.
  
  Given $\varphi\in\Sigma(X)$, let $K_\varphi\subseteq \H(X)$ be the
  sublattice\footnote{In this note, I take lattices to be what are
    sometimes called ``bounded lattices'', in other words posets that
    have all finite meets and finite joins and hence in particular
    (since the empty set is finite) have least and greatest
    elements.}  generated by the finitely many elements of the form
  $[\psi]$ where $\psi$ is a sub-term of the term $\varphi$ (write
  $\psi\in \varphi$ to indicate this). Since finitely generated
  distributive lattices are finite, $K_\varphi$ is finite and hence is
  a complete Heyting algebra (although not necessarily a sub-Heyting
  algebra of $\H(X)$). Let $f_\varphi : \H(X) \morphism K_\varphi$ be
  the Heyting algebra morphism induced by the function
  $X \morphism K_\varphi$ that maps $x\in X$ to $[x]$ if $x$ occurs in
  $\varphi$ and to $\bot$ otherwise. It suffices to prove
  \begin{equation}
    \label{eq:60}
    \forall \psi\in\varphi.\; f_\varphi[\psi] = [\psi]
  \end{equation}
  since from \eqref{eq:60} we have that $f_\varphi$ is a surjection
  and ${f_\varphi[\varphi] = \top} \Imp {[\varphi]=\top}$, as
  required.  The proof of \eqref{eq:60} is by induction on the
  structure of $\psi$. The case when $\psi$ is $x\in X$ holds by
  definition of $f_\varphi$. The cases for finite meets and joins are
  immediate, because $K_\varphi$ is a sublattice of $\H(X)$ and
  $f_\varphi$ preserves those operations. Finally, in case
  $(\psi\imp\theta)\in\varphi$, then $f_\varphi$ sends
  $[\psi\imp\theta]$ to the Heyting implication of $[\psi]$ and
  $[\theta]$ in $K_\varphi$, that is, to the greatest 
  $h \in K_\varphi$ with $h \meet [\psi] \leq [\theta]$ in $\H(X)$;
  and since $\psi\imp\theta \in \varphi$, that greatest $h$ is 
  $[\psi\imp\theta]$.
\end{proof}

\subsection{Finitely presented Heyting algebras}
\label{sec:quo}

Recall that quotients of a Heyting algebra $H$ correspond to
\emph{filters}, that is, subsets $F\subseteq H$ that are closed
upwards and closed under finite meets. Given such a filter, the
corresponding congruence on $H$ identifies $h$ and $k$ just when
$h \bimp k \in F$. Write $H/F$ for the quotient Heyting
algebra and $h/F$ for the equivalence class of
$h\in H$ in $H/F$. Given a morphism
$f:H\morphism K$ in $\Heyt$, its \emph{kernel} is the filter 
\begin{equation}
  \label{eq:16}
  \ker f \defeq \{h \in H \mid f(h)=\top\}
\end{equation}
and $f$ has an epi-mono factorization as $H\epi H/\ker f \mono K$.

Finitely generated filters $F$ are \emph{principal}, that is of the
form ${\uparrow_H}(h) \defeq \{k \in H \mid h \leq k \} $ for some $h$
(since if $F$ is generated by a finite subset $S$, we can take
$h = \Meet S$). In this case we write $H/F$ just as $H/h$ and write
the equivalence class $k/F$ just as $k/h$.  The quotient $H\epi H/h$
is isomorphic to the morphism
\begin{equation}
  \label{eq:21}
  q_h : H\morphism \down_H(h)
\end{equation}
where $\down_H(h) \defeq \{h' \in H \mid h' \leq h \}$ and
$q_h$ maps each $h' \in H$ to $h'\conj h$.

Like any category of algebras for a finitary algebraic theory, $\Heyt$
is a locally finitely presentable category~\cite{GabrielP:lokpk}.  The
finitely presented objects are quotients of free objects by finitely
generated congruences; and a congruence generated by equations
$h_1 = k_1,\ldots,h_n=k_n$ corresponds to the principal filter
$\up(\Meet\{h_i \bimp k_i \mid i = 1,\ldots,n\}$). Thus we have that
\emph{finitely presented Heyting algebras} are the ones isomorphic to
$\H(X)/h$ for some finite set $X$ and some $h\in\H(X)$.

\subsection{Residual finiteness}
\label{sec:resf}

Given a Heyting algebra $H$, say that $c\in H$ \emph{has finite
  height} in $H$ if $\down_H(c)$ is a finite set. Write $\Fin(H)$ for
the subset of all such $c$.

\begin{proposition}
  \label{prop:1}
  If $f:H\morphism K$ is a morphism in $\Heyt$ from a finitely
  presented Heyting algebra to a finite one, then $\ker f$
  is a principal filter ${\uparrow_H}(h)$ for some $h\in \Fin(H)$.
\end{proposition}
\begin{proof}
  We can assume that $f$ is a surjection, since otherwise we can
  replace $K$ by the image of $f$.  If $\ker f $ is principal, say
  $\ker f = {\uparrow_H}(h)$, then (as remarked in
  section~\ref{sec:quo}) $\down_H(h)$ is isomorphic to the finite
  Heyting algebra $K$; hence $h$ has finite height. So we just have to
  show that $\ker f$ is a principal filter.
  We can also assume $H$ is a finitely generated free algebra, that
  is, $H = \H(X)$ for some finite set $X$; because given a surjective
  morphism $q:\H(X)\epi H$, if $\ker(f\comp q)=\up_{\H(X)}(k)$ then
  $\ker(f)= \up_H(q(k))$.

  Given that $f$ is a surjective function $\H(X)\epi K$, for each $k$
  in the finite set $K$ we can find some $\varphi_k \in \Sigma(X)$
  with $f[\varphi_k] = k$. Consider the following finite collection of
  algebraic terms in $\Sigma(X)$:
  \[
    x \bimp \varphi_{f[x]} \qquad
    \varphi_{\Meet F} \bimp \Meet\{\varphi_k \mid k\in F \} \qquad
    \varphi_{\Join F} \bimp \Join\{\varphi_k \mid k\in F \} \qquad
    \varphi_{k\imp k'} \bimp \varphi_k \imp \varphi_{k'}
  \]
  where $x$ ranges over the finite set $X$, $F$ ranges over finite
  subsets of $K$, and $k$ and $k'$ range over $K$. Let
  $\varphi_f\in \Sigma(X)$ be the conjunction term for this finite set
  of terms. Since $f$ is a morphism of Heyting algebras, by
  construction of $\varphi_f$ we have $f[\varphi_f] = \top$, so that
  ${\uparrow_H}[\varphi_f] \subseteq \ker f$. Furthermore, one can
  prove by induction on the structure of the term $\psi \in \Sigma(X)$
  that
  \[
    \forall \psi \in \Sigma(X).\;\varphi_f \leq (\varphi_{f[\psi]} \bimp
    \psi)
  \]
  Hence if $[\psi]\in\ker f$ then $\varphi_{f[\psi]} = \varphi_\top$,
  so $\varphi_f \leq \varphi_\top \bimp \psi$ and hence that
  $\varphi_f\leq \psi$ (because $\varphi_f \leq \varphi_\top$ by
  construction).  Therefore altogether we have that
  $\ker f = {\uparrow_H}[\varphi_f]$ is a principal filter.
\end{proof}

In view of the above proposition the usual notion of residual finiteness
from universal algebra becomes for Heyting algebras the following:

\begin{definition}
  \label{def:resf}
  A Heyting algebra $H$ is \emph{residually finite} if for all
  $h, k \in H$
  \[
  (\forall c\in \Fin(H).\; h \meet c \leq k \meet c) \Imp {h \leq k}
  \]
  or equivalently, if the morphism
   \begin{equation}
    \label{eq:17}
    \begin{aligned}
      H &\morphism \textstyle\prod_{c\in\Fin(H)} \down_{H}(c)\\
      h &\mapsto (h \conj c \mid c \in\Fin(H))
    \end{aligned} 
  \end{equation}
  is injective, or equivalently, if $\top$ is the join of
  $\Fin(H)$ in $H$.
\end{definition}

\begin{proposition}
  \label{prop:resf}
  Every finitely presented Heyting algebra is residually finite.
\end{proposition}
\begin{proof}
  Note that if $H$ is residually finite, then for any $h \in H$
  \[
  h = h \conj \top = h \conj \Join\Fin(H) = \Join\{h \conj c \mid c
  \in \Fin(H)\} = \Join \{c \in \Fin(H) \mid c \leq h\} =
  \Join\Fin(\down_H(h))
  \]
  so that $\down_H(h) \iso H/h$ is also residually finite. So it suffices to
  prove the proposition when $H=\H(X)$ with $X$ finite. In that case,
  that \eqref{eq:17} is injective follows by combining
  Propositions~\ref{prop:fmp} and \ref{prop:1}.
\end{proof}

\begin{corollary}
  \label{cor:inf}
  A finitely presented Heyting algebra that is infinite has infinitely
  many elements of finite height.
\end{corollary}
\begin{proof}
  If $H$ is finitely presented then by Proposition~\ref{prop:resf}, it
  is residually finite. So if $\Fin(H)$ is finite, then the codomain
  of the injective function \eqref{eq:17} is a finite product of finite
  Heyting algebras and hence its domain is also finite.
\end{proof}

\subsection{Coprimes}
\label{sec:cop}

Given a Heyting algebra $H$, an element $c\in H$ is \emph{coprime} if
for all finite $S\subseteq H$ it is the case that
$c\leq\Join S \Imp \exists h \in S.\;c \leq h$ (so that in particular
$c\not= \bot$). We need the following fact about free Heyting algebras.
\begin{proposition}
  \label{prop:cof}
  In a free Heyting algebra, every element
  is a finite join of coprime elements.
\end{proposition}
\begin{proof}
  For finitely generated free Heyting algebras, this is an immediate
  consequence of~\citep[Theorem~3]{UrquhartA:freha}. I prefer a proof
  using Sierpi{\'n}ski cones (defined in the next section) that works
  for arbitrary free Heyting algebras and that I learned from Peter
  Freyd in 1980. It is given in Appendix~\ref{sec:proc}.
\end{proof}

\subsection{Sierpi{\'n}ski cone}
\label{sec:siec}

The \emph{Sierpi{\'n}ski cone} $\scn H$ of a Heyting algebra $H$
(sometimes called the \emph{Freyd cover}, or the \emph{scone} of $H$)
is obtained as a poset by adjoining a new greatest element ($\top$) to
$H$. Write $(\_)^- : H \morphism \scn H$ for the order embedding of
$H$ into $\scn H$ as the non-$\top$ elements (so that in particular
$\top^- < \top$ in $\scn H$).

This construction is a special case of Artin--Wraith
Glueing~\citep{WraithG:atrg} \citep[section A2.1,
Example~2.1.12]{JohnstonePT:skeett}, which was used by Freyd to give,
among other things, category-theoretic proofs of existence and
disjunction properties of intuitionistic
logic~\citep[p~250]{LambekJ:inthoc} (see Remark~\ref{rem:cop}).

We need the following properties of Sierpi{\'n}ski cones.
\begin{enumerate}
  
\item\label{item:A} $\scn H$ is a Heyting algebra, with the
  Heyting implication for non-greatest elements given by
  \begin{equation}
    \label{eq:1}
    h^- \imp k^- =
    \begin{cases}
      \top &\text{if $h\leq k$}\\
      (h\imp k)^- &\text{otherwise}
    \end{cases}
  \end{equation}
  In fact $\scn H$ is a \emph{local} Heyting algebra, in the sense
  that th set of its non-greatest elements has a join, namely
  $\top^-$:
  \begin{equation}
    \label{eq:6}
    \forall u \in \scn H.\; {u < \top} \Bimp {u \leq \top^-} 
  \end{equation}
  and in particular $\top$ is coprime in $\scn H$.
 
\item\label{item:B} The function $\p : \scn H \morphism H$ that maps
  $\top\in\scn H$ to $\top\in H$ and maps $h^-\in\scn H$ to $h\in H$
  is a morphism of Heyting algebras. Furthermore, every isomorphism
  $\scn H \iso \scn H'$ in $\Heyt$ restricts to an isomorphism
  $H\iso H'$ that respects $\p$:
  \[
    \xymatrix{{\scn H} \ar@{}[r]|-{\iso} \ar[d]_{\p} 
      & {\scn H'} \ar[d]^{\p}\\
    H \ar@{}[r]|-{\iso} & H'}
  \]

\item\label{item:C} The function $\d : H \morphism \scn H$ that
  maps $\top\in H$ to $\top \in \scn H$ and maps each $h < \top$ to
  $h^-$ preserves finite meets and Heyting implications. If $\top$ is
  coprime in $H$ then $\d$ also preserves finite joins and hence is a
  morphism of Heyting algebras.
  
\end{enumerate}

\section{Topos-Inadmissible Heyting Algebras}
\label{sec:topiha}

If $\E$ is a (locally small) topos~\cite{JohnstonePT:skeett}, then its
subobject classifier $\Omega_\E$ is a Heyting algebra object in $\E$
and hence the set of its truth values $\E(1,\Omega_\E)$ is a Heyting
algebra in the category $\Set$ of sets and functions (one that is
isomorphic to the poset of subobjects of the terminal object
$1$). Recall from the Introduction that we call a Heyting algebra
\emph{topos-admissible} if it is isomorphic to one of this form and
\emph{topos-inadmissible} otherwise.

\begin{proposition}
  \label{prop:2}
  The class of topos-admissible Heyting algebras is closed under
  taking quotients.
\end{proposition}
\begin{proof}
  If $F$ is a filter on $\E(1,\Omega_\E)$, then
  $\E(1,\Omega_\E)/F$ is isomorphic to $\E'(1, \Omega_{\E'})$
  where $\E'$ is the \emph{filterquotient topos}~\citep[section A2.1,
  Example~2.1.13]{JohnstonePT:skeett} formed from $\E$ and $F$.
\end{proof}

\begin{theorem}
  \label{thm:main}
  All free Heyting algebras with two or more generators are
  topos-inadmissible.\\ 
  \emph{(Note that $\H(0)$ and $\H(1)$, being complete Heyting
    algebras, are topos-admissible, via toposes of sheaves.)}
\end{theorem}
\begin{proof}
  Note that when $|X|\geq 2$, $\H(2)$ is a retract of $\H(X)$ and in
  particular a quotient of it. So the theorem follows by combining
  Proposition~\ref{prop:2} with Propositions~\ref{prop:admissible} and
  \ref{prop:incompleteness} below.
\end{proof}

Suppose $\E$ is a topos and that $G\subseteq\E(1,\Omega_\E)$ is a
finite set of its truth-values generating $\E(1,\Omega_\E)$ as a
Heyting algebra and containing a distinguished element $a$.  Following
\citep{KocsisZ:twoapf,vanGoolS:gentvt}, consider the following
sentence in the internal higher-order
logic~\cite[chapter~D4]{JohnstonePT:skeett} of $\E$:
\begin{equation}
  \label{eq:25}
  \forall x:\Omega.\; x \in \Omega_G \imp (x\disj (x \imp a))
\end{equation}
where $\Omega_G$ is a term in the internal language expressing the
intersection of all subsets of $\Omega_\E$ that contain each $g$ from
the finite set $G$ and are closed under the operations $\bot$,
$\conj$, $\disj$ and $\imp$.  So $\Omega_G$ is a way\footnote{In the
  terminology of \citep{vanGoolS:gentvt}, this is a ``top-down''
  description using impredicative quantification; since we are not
  assuming that toposes have natural number objects, the
  ``bottom-up'', number-indexed iterative construction is not
  available.} of describing the internal Heyting subalgebra of
$\Omega_\E$ generated by the finite set of truth values $G$. The
sentence \eqref{eq:25} defines a truth value that we call
\begin{equation}
  \label{eq:26}
  \km(\E,G,a) : 1 \morphism \Omega_\E
\end{equation}
For any $h\in \E(1,\Omega_\E)$, we always have
\begin{equation}
  \label{eq:28}
  \km(\E,G,a) \leq h \disj (h\imp a)
\end{equation}
because of the way $\km(\E,G,a)$ is defined from
\eqref{eq:25}. (Recall that we are assuming that $G$ generates
$\E(1,\Omega_\E)$, so one can prove by induction on the structure of a
description of $h\in \E(1,\Omega_\E)$ built up from the elements of $G$ using the
Heyting algebra operations, that $h:1\morphism \Omega_\E$ factors
through $\Omega_G\mono \Omega_\E$.)

\begin{definition}
  \label{def:km}
  Given a Heyting algebra $H$ with distinguished element $a\in H$,
  define the subset $\KM(H,a)$ by:
  \begin{equation}
    \label{eq:24}
    \KM(H,a) \defeq \{h \disj(h\imp a) \mid h \in H\} 
  \end{equation}
\end{definition}

\begin{lemma}
  \label{lem:4}
  With $\E$, $G$ and $a$ as above, if $H\defeq\E(1,\Omega_\E)$ is a
  \emph{finite} Heyting algebra then $\km(\E,G,a)$ is equal to the
  finite meet $\Meet\KM(H,a)$. 
\end{lemma}
\begin{proof}
  (See the proof of \citep[Theorem~2]{vanGoolS:gentvt}.)  We have
  $\km(\E,G,a) \leq\Meet\KM(H,a)$ by \eqref{eq:28}.  Conversely, to
  see that $\Meet\KM(H,a)\leq \km(\E,G,a)$, suppose
  $H=\{h_1,\ldots,h_n\}$ and in the internal language of $\E$ consider
  the subset
  $S=\{x : \Omega_\E \mid (x = h_1) \disj \cdots\disj (x =
  h_n)\}$. This is an internal Heyting subalgebra of $\Omega_\E$
  containing $G$ and hence contains $\Omega_G$ by definition of the
  latter as an internal intersection of all such. Arguing in the
  internal higher-order logic: if $\Meet\KM(H,a)$ holds then we have
  that $h_i \disj(h_i \imp a)$ holds for all $i=1,\ldots,n$; so if
  $x\in \Omega_G$, then $x \in S$, so $x = h_i$ for some $i$ and
  therefore $x \disj (x \imp a)$ holds because it is equal to
  $h_i \disj (h_i \imp a)$. So altogether $\Meet\KM(H,a)$ implies
  \eqref{eq:25}, that is $\Meet\KM(H,a)\leq \km(\E,G,a)$.
\end{proof}

\begin{proposition}
  \label{prop:admissible}
  If a finitely presented Heyting algebra $H$ is topos-admissible,
  then for any $a\in H$ the subset $\KM(H,a)$ in \eqref{eq:24} has a
  meet.
\end{proposition}
\begin{proof}
  Suppose $H = \E(1,\Omega_\E)$ and $G\subseteq\E(1,\Omega_\E)$ is a
  finite set generating $H$. Without loss of generality, we can assume
  $a\in G$. We will show that $\km(\E,G,a)$ is the meet of $\KM(H,a)$.
  In view of \eqref{eq:28}, for this it suffices to show that if
  $m\in H$ is some lower bound for $\KM(H,a)$, then $m\leq \km(\E,G,a)$,
  which by Proposition~\ref{prop:resf} is the same as showing
  $m \conj c \leq \km(\E,G,a) \conj c$ for all $c\in\Fin(H)$.

  It is a fundamental property of toposes that the slice $\E/C$ of a
  topos $\E$ by an object $C$ is again a topos and that the pullback
  functor $\E\morphism \E/C$ is logical;
  see~\citep[section~A2.3]{JohnstonePT:skeett}. When $C \mono 1$ is
  the extension of a truth value $c:1\morphism \Omega_\E$, then
  $\E\morphism \E/C$ corresponds to the translation between theories
  in higher-order logic obtained by adding $c$ as an axiom to the
  theory of $\E$ and
  $\E/C(1,\Omega_{\E/C}) \iso \E(1,\Omega_\E)/c = H/c \iso
  \down_H(c)$. Under this isomorphism
  $\km(\E/C,\{g\conj c \mid g \in G\}, a \conj c)$ corresponds to
  $\km(\E,G,a) \conj c$, because $\E\morphism \E/C$ preserves
  higher-order logic.  Note that the truth values in $\E/C$
  corresponding to $g\conj c$ with $g \in G$ generate
  $\E/C(1,\Omega_{\E/C})$, because the latter is a quotient of $H$ and
  that is generated by $G$. If $c\in \Fin(H)$, then
  $\E/C(1,\Omega_{\E/C}) \iso \down_H(c)$ is finite and so by
  Lemma~\ref{lem:4} (and the fact that $h \imp(a \conj c) = h \imp a$
  when $h\leq c$) we have
  \begin{equation}
    \label{eq:29}
    \km(\E,G,a) \conj c = \Meet \KM(\down_H(c), a\conj c)
    = \Meet\{h \disj((h \imp a) \conj c) \mid h \leq c\}
  \end{equation}
  Now suppose $m\in H$ is some lower bound for $\KM(H,a)$. Then for any
  $h\leq c$ we have $m\leq h \disj(h \imp a)$ and hence also
  $m\conj c \leq (h \disj (h \imp a))\conj c = h \disj ((h \imp a)
  \conj c)$; therefore \eqref{eq:29} gives
  $m \conj c \leq \km(\E,G,a) \conj c$ for all $c\in \Fin(H)$, as
  required.
\end{proof}

\begin{proposition}
  \label{prop:incompleteness}
  When $H$ is the free Heyting algebra on two generators, $a\not= b$,
  the subset $\KM(H,a)$ does not have a meet.
\end{proposition}
\begin{proof}
  Write $H_b$ for the subalgebra of $H$ generated by $b$. Thus $H_b$
  is the image of the morphism $\H(1)\morphism H$ induced by mapping
  the unique element $0\in 1$ to $b$; and that is split by the
  morphism $H=\H\{a,b\} \morphism \H(1)$ induced by mapping $a$ to
  $\bot$ and $b$ to $0$. So $H_b$ is an isomorphic copy of
  $\H(1)$. For each $c\in\Fin(H_b)$
  consider the Sierpi{\'n}ski cone $\scn(\down_{H_b}(c))$. Let
  \begin{equation}
    \label{eq:20}
    f_c: H\morphism \scn(\down_{H_b}(c))
  \end{equation}
  be the morphism from the free algebra $H$ that is induced by the
  function mapping $a$ to $c^-$ and $b$ to $(b\conj c)^-$. Since
  $\scn(\down_{H_b}(c))$ is finite, by Proposition~\ref{prop:1} there
  is some $\hat{c} \in H$ with $\ker f_c = \up_H(\hat{c})$.
  
  The family $(\hat{c} \mid c\in\Fin(H_b))$ has the following properties:
  \begin{itemize}
    \item[(i)] $\forall h \in H.\; {f_c(h) = \top} \Bimp {\hat{c} \leq
        h}$\\
      \emph{Proof}: since $\ker f_c = \up_H(\hat{c})$.
  \item[(ii)] $\hat{c}$ is coprime.\\
    \emph{Proof}: by (i), since the greatest element in a
    Sierpi{\'n}ski cone is coprime.
    
  \item[(iii)] $\hat{c} \not\leq a$\\
    \emph{Proof}: by (i), since $f_c(a) = c^- \not= \top$.
    
  \item[(iv)] $\forall h \in H.\; \hat{c} \leq h \disj (h \imp a)$\\
    \emph{Proof}: if $\hat{c}\leq h$ then
    $\hat{c} \leq h \disj (h \imp a)$; and if $\hat{c}\not\leq h$,
    then by (i) $f_c(h) < \top$ and hence from \eqref{eq:6} we have
    $f_c(h) \leq c^-$, so that
    $f_c(h \imp a) = f_c(h)\imp f_c(a) = f_c(h)\imp c^- = \top$ and so
    by (i) again $\hat{c} \leq (h \imp a) \leq h \disj (h \imp a)$.

  \item[(v)] $f_c$ is a surjection.\\
    \emph{Proof}: since $f_c(\top) = \top$, to see that $f_c$ is a
    surjection it suffices to show that
    $S_b \defeq \{h\in H_b \mid \exists h' \in H.\; (h \conj c)^- =
    f_c(h')\}$ is the whole of $H_b$; but this is the case because
    $H_b$ is freely generated by $b$ and one can check that $S_b$ is a
    Heyting subalgebra of $H_b$ containing $b$ (closure of $S_b$ under
    Heyting implications uses \eqref{eq:1} and the fact that when
    $h\conj c \leq k \conj c$ then
    $((h\imp k)\conj c)^- = c^- = f_c(a)$).

  \item[(vi)] The Heyting algebra morphism
    $H_b \subseteq H \stackrel{f_c}{\morphism} \scn(\down_{H_b}(c))
    \stackrel{\p}{\morphism} \down_{H_b}(c)$ (where $\p$ is as in
    \ref{sec:siec}(\ref{item:B})) is equal to the quotient morphism
    $q_c$ from
    \eqref{eq:21}.\\
    \emph{Proof}: since $H_b$ is freely generated by $b$, it suffices
    to check that both morphisms map $b$ to $b\conj c$.
    
  \item[(vii)] ${\hat{c} = \hat{c'}} \Imp {c =
      c'}$\\
    \emph{Proof}: if $\hat{c} = \hat{c'}$, then
    $\ker f_c = \ker f_{c'}$ and hence by (v), $f_c$ and $f_{c'}$ are
    isomorphic quotients of $H$. So using \ref{sec:siec}(\ref{item:B})
    we get
    \[
      \xymatrix @C=0pt {
        & H_b \ar@{^{(}->}[d] &\\
        & H \ar[dl]_{f_c} \ar[dr]^{f_{c'}}&\\
        {\scn(\down_{H_b}(c))}\ar[d]_{\p} &
        \iso &
        {\scn(\down_{H_b}(c'))}\ar[d]^{\p}\\
        {\down_{H_b}(c)} & \iso & {\down_{H_b}(c')}
      }
    \]
    and hence by (vi) $q_c$ and $q_{c'}$ are isomorphic quotients of
    $H_b$. Hence their kernels $\up_{H_b}(c)$ and $\up_{H_b}(c')$ are
    equal and therefore $c = c'$. 
  \end{itemize}
  Now suppose that $\KM(H,a)$ does have a meet in $H$, $m$ say. By
  (iv) we have that $\hat{c}\leq m$ for all $c\in\Fin(H_b)$.  By
  Proposition~\ref{prop:cof} we have $m = \Join D$ with $D$ a finite
  set of coprimes; so by (ii), $\Fin(H_b)$ can be expressed as a
  finite union
  $\textstyle\bigcup_{d\in D}\{ c \in \Fin(H_b) \mid \hat{c}\leq
  d\}$. But $\Fin(H_b)$ is infinite by Corollary~\ref{cor:inf}, since
  $H_b\cong \H(1)$ is
  infinite~\citep{RiegerL:lattbp,NishimuraI:forvip}. Therefore by the
  pigeonhole principle, for some $d\in D$ there exist $c\not= c'$ in
  $\Fin(H_b)$ with $\hat{c}\leq d$ and $\hat{c'} \leq d$. So we have
  $d\not\leq \hat{c}\imp a$, since otherwise
  $\hat{c} = \hat{c}\conj d \leq a$ contradicting (iii); and similarly
  $d\not\leq \hat{c'}\imp a$. But by assumption on $m$, we have
  $d\leq m\leq \hat{c} \disj (\hat{c}\imp a)$ and
  $d\leq m\leq \hat{c'} \disj (\hat{c'}\imp a)$. So since $d$ is
  coprime, we must have $d\leq \hat{c}$ and $d\leq \hat{c'}$. Therefore
  $\hat{c} = d = \hat{c'}$ and hence $c = c'$ by (vii),
  contradiction. Thus $\KM(H,a)$ cannot have a meet.
\end{proof}

\paragraph{Acknowledgement}

As my first serious test of \emph{Claude} [\url{https://claude.ai}], I
fed Claude Opus~5.5 everything in this note up to the statement of
Proposition~\ref{prop:incompleteness} and asked it for a proof using
the method of Sierpi{\'n}ski cones rather than the universal Kripke
model of \citet{BellissimaF:fingfh}.  The overall structure of the
proof given here was the one it found, but I have cleaned up the
details.


\begin{thebibliography}{30}
\expandafter\ifx\csname natexlab\endcsname\relax\def\natexlab#1{#1}\fi
\expandafter\ifx\csname url\endcsname\relax
  \def\url#1{{\tt #1}}\fi

\bibitem[Abbadini et~al.(2026)Abbadini, Almeida, and Arrieta]{AbbadiniM:topefh}
M.~Abbadini, R.~N. Almeida, and I.~Arrieta.
\newblock A topos for {\'e}tale-finite {Heyting} algebras.
\newblock {\em ArXiv e-prints}, arXiv:2606.03861 [math.LO], June 2026.
\newblock URL \url{https://arxiv.org/abs/2606.03861}.

\bibitem[Aczel(1968)]{AczelP:satit}
P.~Aczel.
\newblock Saturated intuitionistic theories.
\newblock In H.~Schmidt, K.~Sch{\"u}tte, and H.~Thiele, editors, {\em
  Contributions to Mathematical Logic}, pages 1--11. North Holland, Amsterdam,
  1968.

\bibitem[Balbes and Dwinger(1975)]{BalbesR:disl}
R.~Balbes and P.~Dwinger.
\newblock {\em Distributive Lattices}.
\newblock University of Missouri Press, 1975.

\bibitem[Bellissima(1986)]{BellissimaF:fingfh}
F.~Bellissima.
\newblock Finitely generated free {Heyting} algebras.
\newblock {\em Journal of Symbolic Logic}, 51\penalty0 (1):\penalty0 152--165,
  1986.

\bibitem[Darni{\`e}re and Junker(2010)]{DarniereL:belcfg}
L.~Darni{\`e}re and M.~Junker.
\newblock On {Bellissima's} construction of the finitely generated free
  {Heyting} algebras, and beyond.
\newblock {\em Arch. Math. Logic}, 49:\penalty0 743--771, 2010.

\bibitem[F{\'e}r{\'e}e and Shilito(2026)]{FereeH:pitims}
H.~F{\'e}r{\'e}e and I.~Shilito.
\newblock {Pitts} and intuitionistic multi-succedent: Uniform interpolation for
  {KM}.
\newblock {\em ArXiv e-prints}, arXiv:2602.16445 [cs.LO], May 2026.
\newblock URL \url{https://arxiv.org/abs/2602.16445}.

\bibitem[Gabriel and Ulmer(1971)]{GabrielP:lokpk}
P.~Gabriel and F.~Ulmer.
\newblock {\em Lokal Pr{\"a}sentierbare Kategorien}, volume 221 of {\em Lecture
  Notes in Mathematics}.
\newblock Springer-Verlag, 1971.

\bibitem[Ghilardi(1992)]{GhilardiS:frehab}
S.~Ghilardi.
\newblock Free {Heyting} algebras as bi-{Heyting} algebras.
\newblock {\em C. R. Math. Rep. Acad. Sci. Canada}, XIV\penalty0 (6):\penalty0
  240--244, d{\'ecembre} 1992.

\bibitem[Hyland et~al.(1980)Hyland, Johnstone, and Pitts]{PittsAM:trit}
J.~M.~E. Hyland, P.~T. Johnstone, and A.~M. Pitts.
\newblock Tripos theory.
\newblock {\em Math.\ Proc.\ Cambridge Philos.\ Soc.}, 88:\penalty0 205--232,
  1980.

\bibitem[Johnstone(1977)]{JohnstonePT:topt}
P.~T. Johnstone.
\newblock {\em Topos Theory}.
\newblock Number~10 in LMS Mathematical Monographs. Academic Press, London,
  1977.

\bibitem[Johnstone(2002)]{JohnstonePT:skeett}
P.~T. Johnstone.
\newblock {\em Sketches of an Elephant, A Topos Theory Compendium, Volumes 1
  and 2}.
\newblock Number 43--44 in Oxford Logic Guides. Oxford University Press, 2002.

\bibitem[Johnstone(2009)]{JohnstonePT:hocmt}
P.~T. Johnstone.
\newblock Hochas and minimal toposes, 2009.
\newblock URL \url{https://rmi.tsu.ge/tolo3/Presentations/Johnstone_1.pdf}.
\newblock Invited talk at Category Theory 2009, Cape Town, South Africa.

\bibitem[Kocsis(2026)]{KocsisZ:twoapf}
Z.~Kocsis.
\newblock Two applications of the point-free coderivative.
\newblock {\em ArXiv e-prints}, arXiv:2609.29436 [math.LO], September 2026.
\newblock URL \url{https://arxiv.org/abs/2609.29436}.

\bibitem[Lambek and Scott(1986)]{LambekJ:inthoc}
J.~Lambek and P.~J. Scott.
\newblock {\em Introduction to Higher Order Categorical Logic}.
\newblock Cambridge University Press, 1986.

\bibitem[Lawvere(1970)]{LawvereFW:equhcs}
F.~W. Lawvere.
\newblock Equality in hyperdoctrines and the comprehension schema as an adjoint
  functor.
\newblock In A.~Heller, editor, {\em Applications of Categorical Algebra},
  pages 1--14. American Mathematical Society, Providence RI, 1970.

\bibitem[McKinsey and Tarski(1946)]{McKinseyJCC:cloeca}
J.~C.~C. McKinsey and A.~Tarski.
\newblock On closed elements in closure algebras.
\newblock {\em Annals of Mathematics}, 47\penalty0 (1):\penalty0 122--162,
  January 1946.

\bibitem[McKinsey and Tarski(1948)]{McKinseyJCC:thescl}
J.~C.~C. McKinsey and A.~Tarski.
\newblock Some theorems about the sentential calculi of {L}ewis and {H}eyting.
\newblock {\em Journal of Symbolic Logic}, 13\penalty0 (1):\penalty0 1--15,
  March 1948.

\bibitem[Nishimura(1960)]{NishimuraI:forvip}
I.~Nishimura.
\newblock On formulas of one variable in intuitionistic propositional calculus.
\newblock {\em Journal of Symbolic Logic}, 25\penalty0 (4):\penalty0 327--331,
  1960.

\bibitem[Pataraia(2007)]{PataraiaD:higoqt}
D.~Pataraia.
\newblock Higher order {Pitts} quantifiers and toposes.
\newblock Incomplete typescript, 7 pages, 2007.

\bibitem[Pitts()]{PittsAM:heyaho}
A.~M. Pitts.
\newblock Heyting algebras and higher-order logic.
\newblock URL \url{https://www.cl.cam.ac.uk/~amp12/talks/Heyting2025.pdf}.
\newblock Invited talk at \emph{Heyting Day 2025: Models of Intuitionism and
  Computability, Symposium in Honour of Jaap van Oosten}, KNAW, Amsterdam,
  2025.

\bibitem[Pitts(1992)]{PittsAM:intsoq}
A.~M. Pitts.
\newblock On an interpretation of second order quantification in first order
  intuitionistic propositional logic.
\newblock {\em Journal of Symbolic Logic}, 57:\penalty0 33--52, 1992.

\bibitem[Pitts(2002)]{PittsAM:tritr}
A.~M. Pitts.
\newblock Tripos theory in retrospect.
\newblock {\em Mathematical Structures in Computer Science}, 12:\penalty0
  265--279, 2002.

\bibitem[Rieger(1949)]{RiegerL:lattbp}
L.~Rieger.
\newblock On the lattice theory of {Brouwerian} propositional logic.
\newblock {\em Acta Facultatis Rerum Naturalium Universitatis Carolinae}, 189,
  1949.

\bibitem[Simmons(1982)]{SimmonsH:algvcb}
H.~Simmons.
\newblock An algebraic version of {Cantor-Bendixson} analysis.
\newblock In B.~Banaschewski, editor, {\em Categorical Aspects of Topology and
  Analysis}, volume 915 of {\em Lecture Notes in Mathematics}, pages 310--323,
  Berlin, Heidelberg, 1982. Springer.
\newblock ISBN 978-3-540-39041-1.

\bibitem[Urquhart(1973)]{UrquhartA:freha}
A.~Urquhart.
\newblock Free {H}eyting algebras.
\newblock {\em Algebra Universalis}, 3:\penalty0 94--97, December 1973.

\bibitem[{van Gool}(2026{\natexlab{a}})]{vanGoolS:gentvt}
S.~{van Gool}.
\newblock Generation and truth values in a topos.
\newblock Zenodo, October 2026{\natexlab{a}}.
\newblock URL \url{https://zenodo.org/records/23185917}.

\bibitem[{van Gool}(2026{\natexlab{b}})]{vanGoolS:unii}
S.~{van Gool}.
\newblock Uniform interpolation.
\newblock In B.~{ten Cate}, J.~C. Jung, P.~Koopmann, C.~Wernhard, and
  F.~Wolter, editors, {\em Theory and Applications of Craig Interpolation},
  pages 303--330. Ubiquity Press, 2026{\natexlab{b}}.
\newblock URL \url{https://arxiv.org/abs/2512.15391}.

\bibitem[Wraith(1974)]{WraithG:atrg}
G.~Wraith.
\newblock Artin glueing.
\newblock {\em Journal of Pure and Applied Algebra}, 4:\penalty0 345--348,
  1974.

\bibitem[Xu and Ye(2026)]{XuY:faihot}
Y.~Xu and L.~Ye.
\newblock Failure of higher-order truth within intuitionistic propositional
  logic.
\newblock {\em ArXiv e-prints}, arXiv:2608.26874 [math.CT], 2026.
\newblock URL \url{https://arxiv.org/abs/2608.26874}.

\bibitem[Ye(2026)]{YeL:coecph}
L.~Ye.
\newblock Coexact completion of profinite {Heyting} algebras and uniform
  interpolation.
\newblock {\em ArXiv e-prints}, arXiv:2604.08267v3 [math.LO], 2026.
\newblock URL \url{https://arxiv.org/abs/2604.08267}.

\end{thebibliography}

\appendix

\section{Proof of Proposition~\ref{prop:cof} using Sierpi{\'n}ski
  Cones} 
\label{sec:proc}

Suppose $H =\H(X)$ is a free Heyting algebra. We aim to show that
every element of $H$ is a finite join of coprime elements.

\begin{definition}
  \label{def:size}
  Define the size of a term $\varphi\in\Sigma(X)$ to be the number of
  symbols it contains. Thus
  $\size:\Sigma(X)\morphism \NNO$ satisfies:
  \begin{equation}
    \label{eq:15}
    \begin{aligned}
      \size(x)
      &= 1 &&(x\in X)\\
      \size(\Join S)
      &= 1 + \textstyle\sum_{\varphi \in S} \size(\varphi)\\
      \size(\Meet S)
      &= 1 + \textstyle\sum_{\varphi \in S} \size(\varphi)\\
      \size(\varphi\imp\psi)
      &= 1+\size(\varphi)+\size(\psi)
    \end{aligned}
  \end{equation}
  Then define the size of an element $h$ of the free algebra
  $H = \H(X)$ to be the minimal size of a term that represents it:
 \begin{equation}
   \label{eq:10}
   \size(h) \defeq \min\{\size(\varphi) \mid h = [\varphi] \} 
 \end{equation}
 Say that $h$ is \emph{formally $\Join$-reducible} if $h = \Join S$
 for some finite subset $S\subseteq H$ with
 $\forall k\in S.\; \size(k) < \size(h)$ (so that in particular
 $\bot=\Join \emptyset$ is formally $\Join$-reducible). Say that $h$
 is \emph{formally $\Join$-irreducible} if it is not formally
 $\Join$-reducible.
\end{definition}

Clearly if $c$ is a coprime then it is formally
$\Join$-irreducible. Also, by well-founded induction on the size of
$h$, we have that every $h$ is a finite join of elements that are
formally $\Join$-irreducible. So to complete the proof of
Proposition~\ref{prop:cof}, we just have to prove that every formally
$\join$-irreducible element is a coprime. One can do so using
Sierpi{\'n}ski cones.

Given $c\in H$, consider the Sierpi{\'n}ski cone $\scn(\down_H(c))$
and let $\hat{q}_c: H \morphism \scn(\down_H(c))$ be the Heyting
algebra morphism induced by the following function
$X\morphism \scn(\down_H(c))$:
\[
x\in X \mapsto
\begin{cases}
  \top &\text{if $c\leq x$ in $H$}\\
  (x\conj c)^- &\text{otherwise}
\end{cases}
\]
(where we just write $x$ for $\eta_X(x)\in \H(X) = H$).
Note that composing $\hat{q}_c$ with $\p: \scn(\down_H(c)) \morphism
\down_H(c)$ we get a morphism of Heyting algebras that agrees with
$q_c$ from \eqref{eq:21} on $X$ and therefore
\begin{equation}
  \label{eq:32}
  \p\comp\hat{q}_c = q_c
\end{equation}
Define a binary relation ${\aczel} \subseteq H\times H$ by:
\begin{equation}
  \label{eq:31}
  {c \aczel h} \;\defeq\; {\hat{q}_c(h) = \top}
\end{equation}

\begin{lemma}
  \label{lem:a1}
  An element $c \in H$ is coprime iff $c \aczel c$.
\end{lemma}
\begin{proof}
  When $c$ is coprime, the greatest element of the Heyting algebra
  $\down_H(c)$ is coprime and hence, as noted in
  \ref{sec:siec}(\ref{item:C}),
  $\d : \down_H(c) \morphism \scn(\down_H(c))$ is a morphism of
  Heyting algebras. Its composition with $q_c:H\morphism \down_H(c)$
  agrees with $\hat{q}_c$ on $X$ and hence $\d\comp q_c =
  \hat{q}_c$. Therefore $\hat{q}_c(c) = \d(q_c(c)) = \d(\top) = \top$,
  that is $c\aczel c$.

  Conversely, suppose $c \aczel c$, that is $\hat{q}_c(c)= \top$. If
  $S$ is a finite subset of $H$ with $c \leq \Join S$, then
  $\top = \hat{q}_c(c) \leq \hat{q}_c(\Join S) = \Join\{\hat{q}_c(h) \mid
  h \in S\}$. But as noted in \ref{sec:siec}(\ref{item:A}), $\top$ is coprime in
  $\scn(\down_H(c))$ so that for some $h\in S$ we have
  $\top = \hat{q}_c(h)$. Therefore by \eqref{eq:32},
  $c \conj h = q_c(h) = \p(\hat{q}_c(h)) = \p(\top) = c$, that is
  $c\leq h$. So $c$ is coprime.
\end{proof}

\begin{remark}
  \label{rem:cop}
  Since $\hat{q}_\top$ preserves $\top$ we have $\top\aczel\top$ and thus
  by the Lemma, \emph{the greatest element of a free Heyting algebra
    is coprime}. This is an algebraic way of proving the disjunction
  property of intuitionistic propositional logic.
\end{remark}

\begin{lemma}
  \label{lem:a2}
  The relation \eqref{eq:31} has the
  following properties for all $c,h,k\in H$,
  $x\in X$ and finite subsets $S\subseteq H$:
  \begin{gather}
    {c\aczel h} \Imp c\leq h\label{eq:33}\\
    {c\aczel x} \Bimp {c \leq x}\label{eq:3}\\
    {c\aczel \Meet S} \Bimp {\forall h \in S.\; c \aczel h}
    \label{eq:4}\\
    {c\aczel \Join S} \Bimp {\exists h \in S.\; c \aczel h}
    \label{eq:5}\\
    {c \aczel (h\imp k)} \Bimp
    {c\leq (h\imp k)} \Conj ({c\aczel h} \Imp {c\aczel
      k}) \label{eq:8} 
  \end{gather}
  (Properties \eqref{eq:3}--\eqref{eq:8} characterise the \emph{Aczel
    slash relation~\citep{AczelP:satit}}.)
\end{lemma}
\begin{proof}
  Property \eqref{eq:33} follows immediately from \eqref{eq:32}.
  Property \eqref{eq:3} follows by definition of $\hat{q}_c$ on
  $X$. Properties \eqref{eq:4} and \eqref{eq:5} follow from the fact
  that $\hat{q}_c$ preserves finite meets and joins and that $\top$ is
  coprime in $\scn(\down_H(c))$.  Finally, for \eqref{eq:8} note that
  the left-to-right implication follows from \eqref{eq:33} and the
  fact that $\hat{q}_c$ preserves Heyting implication. For the
  right-to-left implication, suppose $c\leq(h\imp k)$ and
  \begin{equation}
    \label{eq:9}
    {c\aczel h} \Imp {c\aczel k}
  \end{equation}
  To prove $\top = \hat{q}_c(h \imp k)$, since $\hat{q}_c$ preserves
  Heyting implication, it suffices to show that
  $\hat{q}_c(h)\leq \hat{q}_c(k)$. If $\hat{q}_c(k)=\top$ there is
  nothing to prove, so suppose $\hat{q}_c(k) = (k')^-$ for some
  $k' \leq c$ and thus with $k' = k \conj c$ by \eqref{eq:32}.  If
  $\hat{q}_c(h)=\top$ then we can use \eqref{eq:9}; so suppose
  $\hat{q}_c(h)=(h')^-$ for some $h'\leq c$ and thus with
  $h' = h\conj c$ by \eqref{eq:32}. But $c\leq(h\imp k)$ and therefore
  $h' = h\conj c \leq k \conj c = k'$. Hence
  $\hat{q}_c(h)= (h')^- \leq (k')^- = \hat{q}_c(k)$, as required.
\end{proof}

\begin{proof}[Proof of Proposition~\ref{prop:cof}]
  As noted above, we just need to prove that every formally
  $\join$-irreducible element $c \in H=\H(X)$ is coprime, that is,
  satisfies $c\aczel c$. Towards a contradiction, suppose $c$ satisfies
  $c\not\aczel c$, but is formally $\join$-irreducible. Let
  $\varphi\in\Sigma(X)$ be of minimal size among the terms satisfying
  $[\varphi]= c$. Consider the syntactic form of the term $\varphi$:
  
  We cannot have $\varphi = x$
  for some $x \in X$ by \eqref{eq:3}, because $c\not\aczel c$. 
  
  We cannot have $\varphi=\Join S$ for some finite subset
  $S\subseteq\Sigma(X)$, since in that case for each $\psi \in S$,
  $\size(\psi) < 1 + \sum_{\theta\in S}\size(\theta) = \size(\varphi) =
  \size(c)$, so $\size[\psi]\leq \size(\psi) < \size(c)$ and
  $c = \Join\{[\psi] \mid \psi \in S\}$, contradicting the assumption
  that $c$ is formally $\join$-irreducible.
    
  We cannot have $\varphi = \psi \imp \theta$, since in that case
  $\size(\theta) < 1+\size(\psi)+\size(\theta) = \size(\varphi)$; but
  $c\not\aczel c$ and $c = [\psi]\imp[\theta]$, so from \eqref{eq:8}
  we must have $c\aczel[ \psi]$, so from \eqref{eq:33} $c\leq [\psi]$
  and hence $c = [\theta]$, contradicting the choice of $\varphi$ with
  minimal size representing $c$.

  The only remaining possibility is that $\varphi=\Meet S$ for some
  finite subset $S\subseteq\Sigma(X)$. Since $c\not\aczel c$, by
  \eqref{eq:4} it must be the case that $c\not\aczel[\psi]$ for some
  $\psi \in S$. Since $c=[\Meet S]$ we have $c\leq [\psi]$. Write $S'$
  for $S-\{\psi\}$, and consider the syntactic form of the term
  $\psi$.
  \begin{itemize}
  \item It cannot be $x$ for some $x \in X$ by \eqref{eq:3}, because
    $c\not\aczel[\psi]$.
    
  \item It cannot be $\Meet T$ for some finite subset
    $T\subseteq \Sigma(X)$, for in that case we have
    $c = [\Meet(S'\cup T)]$ with
    $\size(\Meet(S'\cup T)) \leq 1 + \sum_{\psi' \in S'}\size(\psi') +
    \sum_{\gamma\in T}\size(\gamma) < 1 + \sum_{\psi' \in
      S'}\size(\psi') + \size(\psi) = \size(\varphi) = \size(c)$,
    contradicting the choice of $\varphi$ with minimal size
    representing $c$.
    
  \item It cannot be $\Join T$ for some finite subset
    $T\subseteq \Sigma(X)$, for in that case we have
    \[
      c = \Join\{[\Meet(S'\cup\{\gamma\})] \mid \gamma\in T\}
    \]
    with for
    each $\gamma\in T$
    \begin{align*}
      \size [\Meet(S'\cup\{\gamma\})]
      &\leq \size(\Meet(S'\cup\{\gamma\}))\\
      &\leq 1 + \textstyle\sum_{\psi' \in S'}\size(\psi') + \size(\gamma)\\
      &< 1 + \textstyle\sum_{\psi' \in
        S'}\size(\psi') + \size(\psi)\\
      &= 1+\textstyle\sum_{\psi' \in S}\size(\psi')
        = \size(\varphi)
        = \size c
    \end{align*}
    contradicting the formal $\Join$-irreducibility of $c$.
  
  \item It cannot be $\theta\imp \theta'$, for then
    $c\not\aczel[\psi] = [\theta]\imp[\theta']$ and
    $c\leq[\psi] = [\theta]\imp[\theta']$; so by \eqref{eq:8} we must
    have $c\aczel [\theta]$, hence $c \leq [\theta]$ by \eqref{eq:33}
    and so $c=[\Meet(S' \cup\{\theta'\})]$ with
    $\size(\Meet(S' \cup\{\theta'\})) < \size(\Meet S) = \size(c)$,
    contradicting the choice of $\varphi$ with minimal size
    representing $c$.
    
  \end{itemize}
  So in all cases we obtain a contradiction and thus it cannot be that
  $c\not\aczel c$ holds when $c$ is formally $\Join$-irreducible, as
  required.
\end{proof}

\end{document}